# Finding Large Primes

Gavriel Yarmish Brooklyn College   Joshua Yarmish Pace University

Jason Yarmish NYU Tandon School of Engineering

Abstract:  In this paper we present and expand upon procedures for obtaining large d digit prime number to an arbitrary probability.  We use a layered approach. The first step is to limit the pool of random number to exclude numbers that are obviously composite.  We first remove any number ending in 1,3,7 or 9. We then exclude numbers whose digital root is not 3, 6, or 9. This sharply reduces the probability of the random number being composite. We then use the Prime Number Theorem to find the probability that the selected number n is prime and use primality tests to increase the probability to an arbitrarily high degree that n is prime. We apply primality tests including Euler's test based on Fermat Little theorem and the Miller-Rabin test. We computed these conditional probabilities and implemented it using the GNU GMP library.

## 1.     Introduction

In 1978 Rivest, Shamir and Adleman (RSA) [1] created the RSA cryptosystem. This system plays a significant role in securing our information on the internet. The security of this system is based on the difficulty in factoring large numbers, which are created by multiplying two very large prime numbers. As we will see, although factoring large integers is a very difficult problem, finding large primes is relatively easy.  The Prime Number Theorem tells us that if we choose a number x at random the chance that x is prime is about $\sim \frac{1}{log_e x}$ i.e. If we pick $log_e x$ numbers at random we expect about 1 to be prime. But how do we know when this x is prime?

We use a probabilistic approach. We choose a large random number of a particular digit size but exclude numerous classes of numbers that we know are composite. We use the prime number theorem to calculate the probability that a random d-digit size number is prime. We then show how this probability increases when we exclude classes of composites. The first step is to limit the pool of random number to exclude numbers that are obviously composite; we remove any number ending in 1,3,7 or 9. We then exclude numbers whose digital root is not 3, 6, or 9. These steps sharply reduce the probability of the random number being composite. We then apply primality tests to increase the probability to an arbitrarily high degree that n is prime. If the test indicates the random number is likely prime we calculate the increased probability. In section 2 we review the Prime Number Theorem and calculate the base probability. We also adjust this probability calculation assuming we exclude numbers ending in 2, 4, 5, 6 or 8. We again adjust the probabilities by excluding numbers with specific digital roots. In section 3 we review the Fermat Primality Test, Euler Test and Miller-Rabin Test. We also calculate exactly how the probability is affected by application of these tests. Finally, in Section 4 we present our application of these tests using C++ and the GNU GMP library.

## 2. Calculating and increasing probabilities

### 2.1 Probability a d-digit prime: the Prime Number theorem

The Prime Number Theorem [2](PNT) gives an asymptotic approximation for

$$\pi(x) = \text{number of prime numbers } \leq x \quad \text{i.e.}$$

$$\lim_{x \to \infty} \pi(x) \frac{\ln x}{x} = 1 \quad or \quad \lim_{x \to \infty} \pi(x) = \frac{x}{\ln x}$$

Thus the odds that a randomly selected number not exceeding x is a prime can be approximated by

$$\frac{\pi(x)}{x} = \frac{(x/\ln x)}{x} = \frac{1}{\ln x}$$

So the number of 75 digit primes is

$$\pi(10^{75}) - \pi(10^{74}) = number\ of\ primes \leq 10^{75} - number\ of\ primes \leq 10^{74}$$
$$= number\ of\ primes\ in\ interval\ (10^{99}, 10^{100})$$

Let N(x)=number of x digit primes.

N(75) can be approximated by $\frac{10^{75}}{\ln 10^{75}} - \frac{10^{74}}{\ln 10^{74}} = \frac{10^{74}}{\ln 10}\left(\frac{665}{5550}\right) = .052037087 * 10^{74}$

In general $N(k) \sim \frac{10^{k-1}}{\ln 10}\left(\frac{9k-10}{k(k-1)}\right)$

There are more precise estimates of $\pi(x)$, one of which is [3]

$$\frac{x}{\ln x - 1} < \pi(x) < \frac{x}{\ln x - 1.1} \quad for\ x \geq 60184$$

With this formula we get

$$.05233970251 * 10^{74} < \pi(75) < .05237015782 * 10^{74}$$

To find the probability that a generated k digit number is actually a prime we divide N(k) by the total number of k digit numbers: $9 * 10^{k-1}$ thus

$$P(A) = \frac{N(k)}{9 * 10^{k-1}} = \frac{\frac{10^{k-1}}{\ln 10}\left(\frac{9k-10}{k(k-1)}\right)}{9 * 10^{k-1}} = \boxed{\frac{9k-10}{9k(k-1)\ln 10}}$$

Where A is the event that the k digit number selected is a prime.

And for k=75 we get

P(A)=.005781899

## 2.2 Excluding obvious composites: Increasing the prime probability

If we use k digits which end in 1, 3, 7 or 9, we can increase the probability that the k digit number is prime. (any multi-digit number ending with 2, 4, 5, 6 or 8 cannot be prime). The number of k digit numbers having a 1, 3, 7 or 9 in the last digits is $9 * 10^{k-2} * 4 = 36 * 10^{k-2}$ so the probability now becomes

$$P(A) = \frac{N(k)}{36 * 10^{k-2}} = \frac{10}{4}\left(\frac{9k-10}{9k(k-1)\ln 10}\right)$$

So we increase the probability by a factor $\frac{10}{4} = 2.5$. We now have P(A)=.014454748

## 2.3 Digital Roots: Further increasing the prime probability

We can further increase the probability of A by avoiding all k digit numbers where digital root are 3, 6 or 9.

Recall [4] that the digital root of a nonnegative integer n (dr(n)) is a single digit obtained by continually summing the digits until a single digit is obtained. dr(n) can be defined using the floor function $\lfloor x \rfloor$ as $dr(n) = n - 9\left\lfloor\frac{n-1}{9}\right\rfloor$ or in terms of congruences

$$dr(n) = \begin{cases} 0 & if\ n = 0 \\ 9 & if\ n \neq 0, n \equiv 0\ mod\ n\ (n\ is\ a\ multiple\ of\ 9) \\ n\ mod\ 9 & if\ n \not\equiv 0\ mod\ 9 \end{cases}$$

Thus if

$$dr(n) = 3 \Longrightarrow n = 9k + 3 \quad for\ k = 0,1,2, \dots$$
$$dr(n) = 6 \Longrightarrow n = 9k + 6 \quad for\ k = 0,1,2, \dots$$
$$dr(n) = 9 \Longrightarrow n = 9k \quad\quad\ \ for\ k = \quad 1,2, \dots$$

So that if dr(n)= 3, 6 or 9, n is divisible by 3 and is composite. If we eliminate all n whose digital root is 3, 6 or 9 we decrease the k digit pool that we can choose n from by 1/3 thus increasing the P(A) by a factor of 3, thus

P(A)=.043364243

i.e. If we restrict our k digit number choice to only those that end in 1,3,7 or 9 and whose digital root is not 3,6, or 9 we can expect about 1 prime in about 23 attempts or $\frac{1}{23}$. But how do we know whether the selected number n is actually prime? We discuss some primality tests [5] [6] [7] which can determine that a given n is definitely not a prime but can tell us that n is a prime with a very high probability i.e. these tests allow for false positives (n is prime when it is actually not) and no false negatives (n is not prime when it actually is).

## 3. Primality tests

### 3.1 Fermat primality test

Fermat's Little Theorem (FLT): If p is prime and a is an integer not divisible by p then $a^{p-1} \equiv 1\ mod\ p$ (and for all a $a^p \equiv a\ mod\ P$). By the contrapositive if $a^{n-1} \not\equiv 1\ mod\ n$ for some a ($a \neq 0\ mod\ n$) then n is composite. Using the contrapositive of FLT we can prove that a number is composite without actually factoring it.

However if n passes the test i.e. $a^{n-1} \equiv 1\ mod\ n$ it is not a proof that n is prime because the converse of FLT is not necessarily true.

Example 1: $n = 561\ a = 2\ and\ yet\ 2^{560} \equiv 1\ mod\ 561, But\ 561 = 3*11*7\ is\ composite$

Example 2: $n = 341\ a = 2\ and\ yet\ 2^{340} \equiv 1\ mod\ 561, But\ 341 = 11*38\ \ is\ composite$

Any composite number m that passes the test is called a "Fermat Pseudoprime to the base a." [8] Thus 561 and 341 are both Fermat Pseudoprimes to the base 2 . In fact there are many numbers n that are composite such that $a^{n-1} \equiv 1\ mod\ n$ for every a such that gcd(a,n)=1 (i.e. a and n are relatively prime). Such numbers are called Carmichael numbers.

FERMAT PRIMALITY TEST: To test if n is prime or composite

1) Choose a number a so that a is not divisible by n
2) Compute $a^{n-1} \bmod n$
3) If $a^{n-1} \equiv 1 \bmod n$ claim "n is a probable prime"
   Repeat steps 1, 2 and 3 with a different a. The probability of n being a prime number increases with each iteration
4) If $a^{n-1} \not\equiv 1 \bmod n$ output "n is definitely composite."

As we noted above Fermat's test fails for all a when n is a Carmichael number. So what are the chances that a Carmichael number will be randomly chosen to be tested – Not likely

1) There are only 7 Carmichael numbers within the first $10^4$ numbers.
2) There are only 585,355 Carmichael numbers within the first $10^{17}$ numbers.

For a randomly chosen odd integer with 100 or more digits the probability that n is a Carmichael is so small that for practical purposes we can consider it to be zero.

Thus, the Fermat test is a good reliable test provided the test number is not a Carmichael number in which case no matter how many times step 3 is repeated the test will conclude that n is a probably prime.

We look for more primality tests that are an improvement over the Fermat Test.

### 3.2   Euler Test

Some Math Background

Fact 1: If a and b are integers, and a, b >0 there exist integers s and t such that gcd(a,b)=sa+tb
   Please see any intro number theory text

Fact 2: If a,b,c are integers and a,b,c>0 such that gcd(a,b)=1 and a|bc then a|c
   i.e. if p is prime and p|bc then p|b or p|c

Fact 3: If p is prime and $x^2 \equiv 1 \bmod p$ then $x \equiv 1 \bmod p$ or $x \equiv -1 \bmod p$
   i.e. if $x^2 \equiv 1 \bmod p$ and $x \not\equiv \pm 1 \bmod p$ then p is not prime.

Proof of fact 2: by fact 1 gcd(a,b) =1=sa+tb for some s,t. c=csa+ctb ⟹a|csa and a|cbt (same a|bc)⟹ $a|c$.

Proof of fact 3: $x^2 \equiv 1 \bmod p \implies x^2 - 1 \equiv b \bmod p$. $(x+1)(x-1) \equiv 0 \bmod p$ Let b=x+1, c=x-1 $bc \equiv 0 \bmod p \implies p|bc \text{ and }$ by Fact 2  p|b or p|c. p|(x-1) or p|(x+1) then $x \equiv 1 \bmod p$ or $x \equiv -1 \bmod p$

Fact 4 from above FTL: If p is prime and a positive integer $p \nmid a$ then $a^{p-1} \equiv 1 \bmod p$.

We note that $a^{\frac{n-1}{2}}$ is the square root of $a^{n-1}$ $\left(a^{\frac{n-1}{2}}\right)^2 = a^{n-1}$ and if $a^{n-1} \equiv 1 \bmod p$, by fact 3

$a^{\frac{n-1}{2}} \equiv \pm 1 \bmod n$, and if not then n is composite.

We have the Euler Primality test:

1) Choose a number a so that a is not divisible by n
2) Compute $a^{\frac{n-1}{2}}$
3) If $a^{\frac{n-1}{2}} \equiv \pm 1 \bmod n$ claim "n is a probable prime."
   Repeat steps 1, 2 and 3 for several different values of a - the probability of n being a prime number increases with each iteration
4) If $a^{\frac{n-1}{2}} \not\equiv \pm 1 \bmod n$ output "n is definitely composite."

We previously saw that Fermat test failed (found n to be pseudoprime when it was in fact a composite number) for n=561 and n=341.

Using Euler's test:

$$\left(2^{\frac{560}{2}}\right) = 2^{280} = 1 \bmod 561$$

$\left(2^{\frac{340}{2}}\right) = 2^{170} = 1 \bmod 341$ so 341 is an Euler pseudoprime to base 2

But $5^{280} = 67 \bmod 561 \not\equiv \pm 1 \bmod 561$ so 561 is composite

$5^{170} = 56 \bmod 341 \not\equiv \pm 1 \bmod 341$ so 341 is composite

Note that 561 is not an Euler pseudoprime base 5 but it is a Fermat pseudoprime base 2. Also note, 341 is a Euler pseudoprime base 2 but not base 5.

Note

1) If the Fermat test concludes that n is composite (i.e. $a^{n-1} \not\equiv 1 \bmod n$) Euler test also finds that n is composite (i.e. $a^{\frac{n-1}{2}} \not\equiv \pm 1 \bmod n$) for if $a^{\frac{n-1}{2}} \equiv \pm 1 \bmod n$ then $\left(a^{\frac{n-1}{2}}\right)^2 \equiv (\pm 1)^2 \bmod n =$ 1 mod n.

2) If Euler test finds n to be composite ($a^{\frac{n-1}{2}} \not\equiv \pm 1 \bmod n$) Fermat test may still fail (mistakenly find n to be pseudoprime). As an example, as above, 561 is composite by the Euler test using a=5 but using the Fermat test $(5^{280})^2 = 67^2 \bmod 561 = 1 \bmod 561$. So for n=561 using a=5 the Euler test correctly finds n composite where the Fermat test fails (it finds 561 a probable prime).

Fact 4: If n is composite it has at least 4 square roots 1 mod n.

Proof: Consider the simple case when $n = pq$, $(p \neq q)$
$$x^2 \equiv 1 \bmod p \Rightarrow x \equiv \pm 1 \bmod p$$

$$x^2 \equiv 1 \bmod q \Rightarrow x \equiv \pm 1 \bmod q$$

We have 4 systems of equations

$$x \equiv 1 \bmod p, \ x \equiv -1 \bmod p, \ x \equiv 1 \bmod p, \ x \equiv -1 \bmod p$$
$$x \equiv 1 \bmod q, \ x \equiv \ \ 1 \bmod q, \ x \equiv -1 \bmod q, \ x \equiv -1 \bmod q$$

Each of these can be solved mod pq using the Chinese remainder Theorem each yielding the square root of 1. It is clear that if n=p,q,$p_2$ we would get more than 4 square roots of 1.

Example: n=15=5*3

$$x \equiv \pm 1 \bmod 5$$

$$x \equiv \pm 1 \bmod 3$$

So we solve
$$\begin{array}{cccc} x \equiv 1 \bmod 5 & x \equiv 1 \bmod 5 & x \equiv -1 \bmod 5 & x \equiv -1 \bmod 5 \\ x \equiv 1 \bmod 3 & x \equiv -1 \bmod 3 & x \equiv 1 \bmod 3 & x \equiv -1 \bmod 3 \\ x = 1 & x = 11 & x = -11 & x = -1 \end{array}$$

So square root of 1 in mod 15 are 1, -1, 11, -11 hence for n composite
$a^{\frac{n-1}{2}} \equiv k \bmod n$ where $k \not\equiv \pm 1$, can nonetheless be a square root of 1. i.e. $a^{\left(\frac{n-1}{2}\right)^2} \equiv k^2 \bmod n \equiv 1 \bmod n$

In this example: $a^{\frac{n-1}{2}} \equiv 11 \bmod 15$ and $a^{n-1} \equiv 11^2 \bmod 15 \equiv 1 \bmod 15$.

In such a case Fermat test finds n probable primes since $a^{n-1} \equiv 1 \bmod n$, but Euler test finds n composite since $a^{\frac{n-1}{2}} \not\equiv \pm 1 \bmod n$.

Unfortunately there are odd composites such that the Euler test $a^{\frac{n-1}{2}} \equiv \pm 1 \bmod n$ for every a with gcd(a,n) = 1. These components are called absolute Euler Pseudoprimes (1729 and 2465 are 2 examples). There are far fewer absolute Euler Pseudoprimes than there are Carmichael numbers.

### 3.3  Miller-Rabin test (MRtest)

The Miller-Rabin test [9] uses fact 3 more extensively than in the Euler test (Fact 3 teaches that if the square root of 1 is not $\pm 1 \bmod n$ then n is composite).

To test n if prime or composite:

1) Choose a such that $2 < a \leq n - 1$
2) Write $n - 1 = 2^k m$ $(m \text{ is odd}, K \geq 1)$
3) In mod n evaluate $b_0 = a^m, b_1 = (a^m)^2, b_2 = (a^m)^{2^2}, b_3 = (a^m)^{2^3}, \ldots b_k = (a^m)^{2^k} = a^{n-1}$
   Note $b_i = b_{i-1}^2$ $i = 1,2,\ldots,k$ i.e. $b_{i-1}$ is square root of $b_i$
4) Consider the first $b_j \equiv 1 \bmod n$ (if $b_j \not\equiv 1 \bmod n$ for all j then n is composite.)
5) If $b_{j-1} \not\equiv \pm 1 \bmod n$ then n is composite (a is called a witness); otherwise n is "probably prime" and is called a strong pseudoprime.

Fact 5: If n is an odd prime one of the following two conditions must hold

a) $b_0 \equiv 1 \bmod n$ or
b) $b_i \equiv -1 \bmod n$ for some i=0,1,2,…,k

Proof: if a) is true then all $b_i \equiv 1 \bmod n$ since $b_i = b_{i-1}^2$ $i = 1,2,3,\ldots,k$
i.e. $b_k = a^{n-1} \equiv 1 \bmod n$ (Fermat Little Theorem)

If b) is true $b_i \equiv -1 \bmod n$ implies $b_j \equiv 1 \bmod n$ $j = i+1, 2, 3, \ldots, k$ again we get $b_k = a^{n-1} \equiv 1 \bmod n$

If neither conditions a or b are satisfied the n is not prime (by contrapositive of Fact 5).

We have the following test: A number a, where $2 \leq a \leq n - 2$ (if a=n-1 then condition b is satisfied $(n-1)^m \equiv -1$ since m is odd) is a witness – the test indicates that n is composite- if:

$b_0 = a^m \not\equiv \pm 1 \bmod n$ and

$b_i \not\equiv -1 \bmod n$ for all i= 0,1,2,…,k

If n is composite, then any a, $1 \leq a \leq n-1$, not a witness is called a <u>Liar</u> (i.e. test indicates n is prime using that a) . Note that a=1 and a=n-1 are trivial liars. The n associated with liar a is called a strong pseudoprime base a.

The probability is less than $\frac{1}{4}$ that test gives wrong answers where n is composite, [10] and if n is composite the probability that the test will indicate n is prime for each of m different as (denoted by $T^m P$) is less than $\left(\frac{1}{4}\right)^m$ i.e. $\text{P(test gives a false positive)} = \text{P(TP|c)} \leq \frac{1}{4}$ and

$$P(test\ failed\ to\ detect\ a\ single\ witness\ in\ m\ attempts) = P(T^m P|c) \leq \left(\frac{1}{4}\right)^m$$

We want to find the reliability of the results i.e. if test indicates prime what is the chance that n is indeed prime. We use Bayes theorem:

$$P(p|T^m P) = \frac{P(p)P(T^m P|p)}{P(p)P(T^m P|p) + P(c)P(T^m P|p)}$$

Where P(p) is the probability that the selected k digit number is prime.

P(c )= 1-P(p)

$$P(T^m P|p) = 1$$

$$P(T^m P|c) \leq \left(\frac{1}{4}\right)^m$$

So, probability $P(p|T^m P)$ that the selected n, which goes through the test m times (each with a different a) and whose test indicates n is prime each time, is indeed prime is

$$P(p|T^m P) = \frac{P(p)(1)}{P(p)(1) + P(c)P(T^m P|c)} \geq \frac{P(p)(1)}{P(p) + P(c)\left(\frac{1}{4}\right)^m}$$

$$> \frac{1}{1 + \frac{P(c)\left(\frac{1}{4}\right)^m}{P(p)}} > 1 - \frac{P(c)\left(\frac{1}{4}\right)^m}{P(p)} =$$

Recall P(k digit selected number is prime)=P(p)= $\frac{9k-10}{k(k-1)\ln 10}$

So P(k digit selected number is composite)P(c )=1-P(p)

As an example For k=75: P(p)= .043364243; P(c )= .956635757 $\frac{P(c)}{P(p)} = 22.66047407$, and 4 iterations m=4.

We have $P(p|T^mP) \geq 1 - \frac{22.66047407}{4^m} = 1 - .088517477 = .911482523$

We get a better than 91% reliability for just 4 repetitions. Has we used the selection process without restricting the pool of k digit numbers we would have

$P(p) = \frac{N(75)}{9*10^{74}} = .005781899$  $P(c) = .994218101$  $\frac{P(c)}{P(p)} = \frac{.994218101}{.005781899} = 171.9535573$

$$P(p|T^mP) > 1 - \frac{P(c)\left(\frac{1}{4}\right)^m}{P(p)} = 1 - \frac{171.9535573}{4^m} = 1 - .671692583 = .328306417$$

Thus using the test for 4 iterations without pre-restricting the pool of random k digit numbers to exclude obvious composite numbers gives us a 32% confidence that the number selected is prime. When we pre-restrict the pool we obtain a 91% confidence that the number is indeed prime.

In either case we can increase other confidence level by performing more iterations of the test.

In summary, if we analyze the $P(p|T^mP) > 1 - \frac{P(c)}{P(p)4^m} = 1 - A$     where $A = \frac{P(c)}{P(p)4^m}$.

Because our confidence is 1-A, the smaller A the more confident we are that indeed n is prime. There are 2 ways to decrease A:  1) increase P(p)  2) increase m. Both methods can be utilized simultaneously.

## 4.     Experimental Results

We implemented this method in C++ using the GNU GMP library for arbitrarily large numbers. We generated 100 random 75-digit numbers, note that we could have used arbitrarily large numbers but thie would have made it unwieldy  to present in figure 1.

As specified in the method above, the random number generator only generated numbers whose last digit ended with a 1, 3, 7 or 9 since all other are certainly composite numbers. The random number generator also excluded all numbers whose digital root was 3, 6 or 9 as these too, are obviously composite.

We then implemented the Rabin-Miller test to determine for each generated random number n if it is composite or probably prime. For each 75-digit random number we looped *m*=10 times. Each time we randomly chose a value *a* of the primality test.  If a particular *a* was a 'witness' then the 75-digit number is composite and the loop ended. If it was not a witness then the probability of the number being prime rose.

In this experiment, the calculation of the probability is the same as calculated in the last section except that m is 10 instead of 4.

We have $P(p|T^mP) \geq 1 - \frac{22.66047407}{4^{10}} = 1 - .000026 = .9999978$. This means that a number that lasted through the full loop of 10 is more than 99.99% likely to be prime.

In Figure 1 below we list all the 100 generated numbers in column 1 together with the output from our program in column 2.

There were 5 primes out of the 100 numbers. We manually checked each of the 100 numbers utilizing an online Prime number checker [11]. As would be expected from the probability, they were all correct.

| Random numbers ending in  1, 3, 7 or 9 not having a prime root of 3, 6 or 9 | P/C |
|---|---|
| 2201117550890485208641986978821781904013607879519698867103820772303619836 19 | COMPOSITE |
| 3165145730262073866310174993400367797972868069174527391342759938635010478 89 | COMPOSITE |
| 2976506427796614235232264497741030700531353468213655501344670189638094533 83 | COMPOSITE |
| 6169338965929670297667096415624183456364452441553247912915732178012021069 73 | COMPOSITE |
| 9098541948144114279602577348033333435682437045400550377170049127594542990 93 | COMPOSITE |
| 2287283480405884383226616295378697856680834597031617486956047953739572303 29 | COMPOSITE |
| 6948916299172506350893942613582518451331602502383001606035790933729560479 79 | COMPOSITE |
| 7983899401325653265219748981058848656872605551061526144654799305417324157 49 | COMPOSITE |
| 7346623945209772160424979255916750208980698856421623061845702573677900155 11 | COMPOSITE |
| 7926946753813510871064819849604511504242970414207893316374940052090979768 11 | COMPOSITE |
| 3629291800854444352350041188521181511278286755240974388282648601372996526 81 | COMPOSITE |
| 5307703592507184388142464222782748809786503392873745908269752872949912257 01 | PRIME |
| 2288411804551894023770611669351526472297141495298344859974560917357747274 39 | COMPOSITE |
| 4294377971459136480898604711446300537961224968687300245241662647224407100 61 | COMPOSITE |
| 8536093523038025447223401245419248509395970022143640919781272215385498776 09 | COMPOSITE |
| 3589978818563389345420392466462492586047112188574211453968382782381524467 77 | COMPOSITE |
| 1263894962465110759627087802831056186799411057784203955472899953096991763 39 | COMPOSITE |
| 9710373700223717121341154629472739617405746439577168825514135488071129696 27 | COMPOSITE |
| 4247088576177462775775303329797575509320113858571155297293020092830392654 09 | COMPOSITE |
| 1621377564594033437800090748171994413790923333173305760029234221871719520 93 | COMPOSITE |
| 4442782656595684423161368340503179421168342644336882440452093650253542680 09 | COMPOSITE |
| 6639708963138376750386444669949925219833971521025511356779898484567937371 37 | COMPOSITE |
| 4568637793414178163808806207939255106809499677256792887829197090382418941 43 | COMPOSITE |
| 1220912408818501772462360005115986420314178114971438614296332566000495825 59 | COMPOSITE |
| 6241897463358961026957226283632017856656530863499427250475314507741293529 61 | COMPOSITE |
| 7712817058408890159271771692559167981924008874976911457828937981578130768 03 | COMPOSITE |
| 8976565954700505589104879145237940681230583428438489584934035074387742072 49 | COMPOSITE |
| 8139604498267749665528597919325289938840819894313128715119865232736641701 91 | COMPOSITE |
| 5814577741898074267957999342164431004096043643538742962600154722915019141 29 | COMPOSITE |
| 7080337099862642019875607484009678247188500835889853919171416120365029048 09 | COMPOSITE |
| 5198805812887543278770841025398656508077323938134389460017596509178643074 53 | COMPOSITE |
| 3436762434675400201200252664495448584004584421641237063579282891548251109 23 | COMPOSITE |
| 1953182112747829324819772335819569487800196229591688627995677405537806636 53 | COMPOSITE |
| 3670950286374787012751599912476024209392195128098499387254253481298114695 9 | COMPOSITE |
| 2347734513754800420233428766928830188833442879096695412156588824069491536 01 | COMPOSITE |
| 8207776440483898378027552404774690591725960194832606718319153050722487254 87 | COMPOSITE |
| 1731783604201757419227068060911433755591686706889835735379534992745274569 41 | COMPOSITE |
| 352587286309648925620598712418753401635877843845564757947194151752536122989 | COMPOSITE |
| 412950554901733866067791437975584987392378735868105895460114886561435700071 | COMPOSITE |
| 959886981933527743690278821104811050770316395823045945694995200839978888717 | PRIME |
| 341605811597851734714184079671083178486352952888206177718090132624515116173 | COMPOSITE |
| 307456279851799419993320572576491577860097603001423781683625128262252810113 | COMPOSITE |
| 861561897381395770893135910500903122705574134773006139426043954898444629991 | COMPOSITE |
| 503043474954418659072834784587266464030923885408262316560401342104239802503 | COMPOSITE |
| 385488806759210977918039435293370004656011701414600351639467339195422057487 | COMPOSITE |
| 722631310385150056748134806905156009811034649771389351939398298743725295183 | COMPOSITE |
| 479311386137473057847145446762229293837435424407139909819876222843994090211 | COMPOSITE |
| 546054883813632519802129209119000112900797760300882087528846621278988335211 | COMPOSITE |
| 415525630693953428877535748401640271251635046197517199604055067073222508219 | COMPOSITE |
| 603125758107624397407654397237491049083387865719126764580414303576596713083 | COMPOSITE |
| 583539337449809254289135866339767797925150593991581110371969691979886323451 | COMPOSITE |
| 970554048765379950099878942382471737856477571005572737458050155240809678751 | COMPOSITE |
| 348227051709814021909041956447390667356328526078675103443936571609856871037 | COMPOSITE |
| 694950471024255134687612976359870348842798597802845559630435991293303392 61 | COMPOSITE |
| 157423491599424468734023346843188611724049090342403385286658888878304272399 | COMPOSITE |
| 917584350898346548463596836624257391824652266352264227659601996598234834117 | COMPOSITE |
| 895259462251963476661405129252500043295719390189716042074160601038502323949 | COMPOSITE |
| 863166595859187539559487597693695481462549734437685776238808198078769781587 | PRIME |
| 526184843831314008459333160200311499487172853871017555184649142723492546481 | COMPOSITE |
| 479326723304565047177328527955518109007318814704352660333248548138288931679 | PRIME |
| 110008376632033067344845193045657169010496567129211768638523876688977273137 | COMPOSITE |
| 381440141375024437814726453792125523269301978174116289526117982114730737551 | COMPOSITE |
| 422948321072854225282037695177173016516830019910100646207782969475815987133 | COMPOSITE |
| 838620866191138137947964921770866680707319753439519956032719941302108356233 | COMPOSITE |
| 107945301674456929336803972681175285351771205506654310730839720485578389159 | COMPOSITE |
| 156223204426790826447373761311060664605773709031508698457458616230144501837 | COMPOSITE |
| 400179301730392495700641494600691033741713485289781916132583440537330722833 | COMPOSITE |
| 915657954086756421559468530071783614938793609798076500495539430064737778003 | COMPOSITE |
| 594448983657375633169165352649593609931155382421533666918031046317067730723 | COMPOSITE |
| 477692672274248075708098983914276440423141053172459031966663065076329113929 | COMPOSITE |
| 949769292271604294351239995503945174192472297405684337451005446692314757619 | COMPOSITE |
| 613847708008962578802649202104141078126675871742194122186900610339612648727 | COMPOSITE |
| 709266160958784583098609638436520716427007259075254119727564160022932021703 | COMPOSITE |
| 230135614018136371236965876307158008414058093569443218891285514814338708479 | COMPOSITE |
| 657478890218574632792043905690832416428358688024319385469727440435740882191 | COMPOSITE |
| 157300711400804888314520708353761570038030583655985581102191771567300849 19 | COMPOSITE |
| 779741615978436021751466346948592901426306315350136574623965824494526375187 | COMPOSITE |
| 701789680806724462358318295026552509026028326594548129827589179467659839267 | COMPOSITE |
| 412824563462620434268098535380828449545125633497015080418609199682933667279 | COMPOSITE |
| 569876037048020585710173924183063843680721382061262556345628090001572699947 | COMPOSITE |
| 977107084071477189293704278613322644403584591133499920044751704893231936463 | COMPOSITE |
| 679315585123345496444608369727813396232326550724527262185260118469592705007 | COMPOSITE |
| 933125324909126920431844614409346071377046242721562254728535823240548425027 | COMPOSITE |
| 965229968255818197619948058469725477602160813673368970795272160032015649031 | COMPOSITE |
| 858542030825079594829661820038976573273390740452998345407830062790692449747 | COMPOSITE |
| 321279061857594270472380594849921324045717746083616240909682967183698020513 | COMPOSITE |
| 236641137703197634404249913390784297359237806907513531411853109806073510991 | COMPOSITE |
| 533165600339804012311571604245660003767767556755599502252896524658865752217 | COMPOSITE |
| 732569361386749771688976227915195141816417532500665956107575227442657865983 | COMPOSITE |
| 795427619760864700069797568013984567474777365311527984128398681881223645973 | COMPOSITE |
| 563959233538003537991201091930755280001354376741888300005859867474093013703 | PRIME |
| 647818539577676979456317737428102604595923089788658206175891170241705123314 | COMPOSITE |
| 390401201810491630562663741867686158210481604862027789019034578889017406527 | COMPOSITE |
| 733424032914570080566935736570334318279299263370756797068964931871382996361 | COMPOSITE |
| 219698239961023941110824254512655034618396395829783400460815104116018281651 | COMPOSITE |
| 874465361286900826168956909053063174845040354416579211732236350022372120191 | COMPOSITE |
| 574263285577390292734496054058241003577278225607947522871270374798271731527 | COMPOSITE |
| 337378475285907785332816956280622555263570773284065327930232439542520335553 | COMPOSITE |
| 923858854417266802249920620513613034217951635456501784026184652148649095851 | COMPOSITE |
| 120314323055938180630442814340198936645018097929076775144963383854569008811 | COMPOSITE |

**Figure 1: 100 75-digit random numbers labeled whether they are prime of composite**